%% file: ALFS.tex
\documentclass[12pt, a4paper]{article}

\usepackage{amsmath}\multlinegap=0pt
\usepackage[applemac]{inputenc}
\usepackage{amsthm}
\usepackage{amssymb}
\usepackage[francais,english]{babel}
\usepackage{pdfsync,color}
\usepackage{subfigure}
\usepackage{graphicx}
\usepackage{graphicx,pict2e}

\numberwithin{equation}{section}
\theoremstyle{plain}
\newtheorem{thm}{Theorem}[section]
\newtheorem{coro}[thm]{Corollary}
\newtheorem{prop}[thm]{Proposition}

\theoremstyle{definition}

\theoremstyle{remark}
\newtheorem{rem}[thm]{Remark}

\newcommand{\ve}{\varepsilon}
\newcommand{\ev}{\mathbf{e}}

\newcommand{\R}{\mathbb{R}}

\newcommand\pref[1]{(\ref{#1})}

\let \eps\varepsilon

\def\<#1,#2>{\left<#1,#2\right>}
\let\bar\overline

\newcommand\Om{{\Omega}}

\newcommand\PD{{(\mathcal{P}_{\mathcal{D}})}}

\newcommand\ovv{{\overline v}}
\newcommand\ovw{{\overline w}}

\newcommand\tet{\theta}
\newcommand\ud{\frac{1}{2}}
\newcommand\al{\alpha}
\newcommand\K{\mathcal{K}}
\newcommand\PR{{(\mathcal{P}_{\mathcal{R}})}}

\title{On the strategic use of risk and undesirable goods in multidimensional screening}
\author {A. Lachapelle, F. Santambrogio\thanks{\scriptsize CEREMADE, UMR CNRS 7534, Universit\'e Paris IX Dauphine, Pl. de Lattre de Tassigny, 75775 Paris Cedex 16, FRANCE
\texttt{lachapelle@ceremade.dauphine.fr}, \texttt{filippo@ceremade.dauphine.fr}}}

\begin{document}

\maketitle

\begin{abstract}
A monopolist sells goods with possibly a characteristic consumers dislike (for instance, he sells random goods to risk averse agents), which does not affect the production costs. We investigate the question whether using undesirable goods is profitable to the seller. We prove that in general this may be the case, depending on the correlation between agents types and aversion. This is due to screening effects that outperform this aversion. We analyze, in a continuous framework, both 1D and multidimensional cases.

\end{abstract}

\bigskip
\noindent
{\bf Keywords:} principal-agent problem, adverse selection, calculus of variations, nonlinear pricing, random contracts.

\section{Introduction}\label{intro}

Principal-Agent problems and the theory of multidimensional nonlinear pricing (see for instance \cite{BolDew} for an overview of the topic)  have known rapid advances at the end of the XX$^{th}$ century with the works of Armstrong \cite{Arm} and Rochet \& Chon\'e \cite{rc}.  This framework considers  a partially informed monopolist (or principal), selling multiproducts to a continuum of agents having multidimensional types. This is the classical adverse selection case of informational asymmetries, involving multidimensional screening. The starting point of the paper, and its motivation, is the following very natural and fundamental question. Suppose that a risk neutral monopolist may ``randomize'' the goods, facing risk averse agents (so that random contracts are undesirable from their viewpoint), then what will happen in the end? Should the seller actually offer randomized goods (i.e. lotteries) or not? Will adverse selection force him to produce those goods preferred by the consumers, or will it be convenient to use this extra possibility in order to make profit?

Although these questions are not novel, only few answers can be found in the literature. A simple case where risk aversion is introduced in a scalar problem involving two agents can be found in \cite{mr_1} or \cite{BolDew}. In such cases, it is shown that the monopolist benefits by randomizing contracts as soon as the difference between both agent's aversion is large enough. Another work on $n$ risk averse agents facing a not informed monopolist has been developed by Maskin \& Riley in \cite{mr_2}, where they focus on auctions on a single and indivisible good. We could also mention the work of Arnott \& Stiglitz (see \cite{as}) for a comparison between random and deterministic  insurance contracts that a monopolist offers - again - to a finite number of buyers with unknown preferences. More recently, and in a more general framework, Carlier proved in \cite{gcthese} an existence result of optimal stochastic contracts in a relaxed version of what we call the Rochet-Chon\'e problem (i.e. the multidimensional screening problem). However, the author does not discuss the profitability of introducing lotteries.

It seems, according to this literature, that this question has been analyzed so far only in discrete models (a finite number of different agents), and often with zero production costs. The novelty of our paper is that we find the same kind of answers in a continuous framework {\it à la} Rochet \& Choné, thanks to a calculus of variations formulation.

Before entering the mathematical details of the principal's optimization problem, we discuss some assumptions in a rather qualitative way. The agents will be supposed to be described by a type $\tet$ and the goods the monopolist may produce by a quality $x$. Both $\tet$ and $x$ belong to a same Euclidean space $\R^d$ (or $\R_+^d$) and the utility that an agent of type $\tet$ derives from consuming a good $x$ is supposed to be the scalar product $\tet\cdot x$. Each coordinate $\tet_i$ can be viewed as a marginal utility to the $i^{th}$ quality of the good. Besides the utility function, the other element of the model is the production cost $C(x)$. This cost is supposed to be convex in $x$.
We consider, as in \cite{rc}, the purely quadratic production cost, i.e. $C(x)=\ud |x|^2$. Both these choices (linear utility and quadratic production costs) are the most classical and the most significant ones, since they already capture all the qualitative features of the solution of more general cases.

We then consider the fact that goods qualities may include a random part, and we will describe it by introducing a new variable $y\in \R_+$ standing for instance for their variance (or for the standard deviation). More precisely, we think of a good as a random variable $X$ where $x=E[X]$ stands for the mean and $y=E[|X-x|^2]=E[X^2]-|E[X]|^2$ is a scalar variance. We suppose that consumer agents have still a linear utility function of the form $\theta\cdot x +\alpha y$, where $\alpha<0$ is a new parameter in the description of agents' type, standing for their risk aversion (see Section \ref{pb} for a more detailed discussion). Notice that we are assuming that the utility is linearly decreasing in $y$, which means that $y$ must be chosen as a parametrization of the randomness of the good such that agents are linearly averse with respect to it. If by chance in modeling a situation we find that agents are linear with respect to the standard deviation instead of the variance, we should call $y$ the standard deviation. But this is, up to now, just a matter of naming variables.

\bigskip
The role played by $y$ would be more important if looking at the production cost. Due to the random components in this extended model, we will call it $C_R(x,y)$. We have in mind at least two toy models for the dependence of $C_R$ on $y$. 

\smallskip
\noindent
$\bullet$ There are some cases where $C_R$ depends increasingly on $y$, say $C_R(x,y)=|x|^2+\lambda y$ for $\lambda>0$, or $|x|^2+\lambda y^2$. This is for instance the case where the monopolist sells black boxes containing an unknown collection of goods, randomly sorted according to a certain distribution (think of closed boxes containing a dozen of wine bottles, randomly selected by the seller). If the cost for each of them is quadratic, the cost for producing one of such boxes is proportional to $E[X^2]$ and not to $|E[X]|^2$. This means that the cost is of the form $|E[X]|^2+\left(E[X^2]-|E[X]|^2\right)=|x|^2+y$.

\smallskip
\noindent
$\bullet$ There are also cases where $C_R$ depends decreasingly on $y$, say $C_R(x,y)=|x|^2+\lambda y$ for $\lambda<0$. This happens when the producer needs special machines for guaranteeing more precise goods qualities, and he does not really chooses to introduce variances, even if he can pay for reducing it. For the model to make sense, the monopolist must be able to manage at the same time, and at his will, the production of possibly different levels of variance, exactly as he needs to be able to produce different levels of quality. This is surely reasonable in large-scale firms, using different machines in different factories.

\bigskip
In this paper we will mainly focus on the case $\lambda=0$, i.e. $C_R(x,y)=|x|^2$. This case can be obtained as a combination of the two previous effects, but more important reasons suggest to consider it as key case.\\
First of all, we already pointed out that we use a parametrization of the variable $y$ so as to guarantee that the utility functions of agents are linear. If $C_R$ really depends on $y$, when specifying a dependence (for instance $C_R(x,y)=|x|^2+\lambda y$) we would be fixing a very specific case where both production costs and utilities are linear in the same parametrization, while there is no modeling reason to suppose that.\\
Second, the case $\lambda>0$ is not a new problem : it is just the Rochet-Choné Problem (see \cite{rc}) with an extra dimension. The only point which was not investigated in \cite{rc} is the fact that we look at consumers in $\R_+^d\times \R_-$ and goods in $\R_+^d\times \R_+$, i.e. we explicitly account for a characteristic that is negatively perceived by agents, but the formalism would be the same.\\
Third, we will see that the main results of the paper consist in examples and conditions where the optimal strategy for the risk neutral principal is exactly to offer goods with nonzero values of $y$ (i.e. lotteries). This makes the case $\lambda<0$ less interesting. Indeed, it is clear that whenever there exist examples where the seller offers costless lotteries, the profit he gets from lotteries increases if $\lambda<0$.

\bigskip
Moreover, it is clear from the choice $\lambda=0$ that the label ``variance'' does not play any special role any more, in the sense that this last characteristic may be thought of as any extra feature of the good that consumers dislike. It is easy to re-think at the whole problem in these terms : if the monopolist may produce goods with qualities $(x,y)$, and the value of $y$ does not affect the production cost, but negatively affects the utilities of the consumers, will he use this ``undesirable'' quality $y$? \\
It is widely accepted that examples of this kind have been used in real life by some companies. It was for instance the case of third-class seats in some trains, which were just more unconfortable than second-class ones without saving space for the railway company. It seems a common strategy of the firms to create low-quality goods just to attract consumers that would otherwise not enter the market.
For these reasons, we will from now on refer to the quality $y$ as the {\it undesirable quality} and to $\al$ as the {\it aversion parameter} (and not specifically risk aversion). And, for the sake of generality, we will drop the random terminology (risk, lotteries, \dots).\\
Actually, thanks to these last considerations, the problem we study and the (partial) answers we give may switch very easily from the random interpretation to the more abstract and general one of an undefined ``undesirable'' quality. This is one of the main features of this very simplified model, which does not want to compete with more sophisticated ones that have been investigated concerning Principal-Agent problems involving random contracts. It is the case of the problems studied by Carlier in \cite{gcthese} (where the author focuses on existence issues for more general utility functions and ``randomness", i.e. not only considering mean and variance of stochastic contracts), as well as the work of Carlier, Ekeland \& Touzi (\cite{CarEkeTou}), in which, nevertheless, the principal cannot control the level of the undesirable quality $y$.

\bigskip
In section \ref{pb} we will first discuss more precisely the {\it maximization problem} for the principal, which we have not yet introduced. We will then jump to the extended problem (i.e. adding the components $y$ and $\al$) and present some considerations on it. Section \ref{1D} and \ref{2D} will be devoted to the description of some classes of examples where the solution of the extended problem really uses $y$ (hence it does not coincide with the solution of the classical problem), in the scalar and the multidimensional cases, respectively. The paper will be closed by a conclusive section with comments on the result we got.

\section{Classical and extended problems}\label{pb}

We start by recalling the classical problem studied by Rochet and Chon\'e, including its powerful calculus of variations formulation. 

\bigskip
\noindent
{\bf General setting} Let us see more in details the very famous problem studied by Mussa \& Rosen (see  \cite{mr}) in the one dimensional case and by Armstrong and Rochet \& Chon\'e for the case of multidimensional types (see respectively \cite{Arm,rc}). It is a nonlinear pricing model in a monopolistic situation with adverse selection. The principal (also the monopolist, the firm, the seller...) knows the global distribution $f$ of agents types on a convex open bounded domain $\Om \subset\R^d_+$. This domain will be supposed for simplicity to be a rectangle and the distribution of agents over it will be described through its density, i.e. a positive function on $\Omega$ with integral equal to $1$, that we still denote by $f$. In this way we consider a continuum of buyers. Any agent with type $\tet \in \Om$  derives the utility $u_{\tet}(x,t):=\tet\cdot x-t$ from consuming a good of quality $x \in Q\subset \R^d_+$ at price $t$. The monopolist has a strictly convex cost $C: x \in Q \rightarrow C(x)$ for producing a good of quality $x$. 

The optimization of the monopolist's strategy is performed on possible {\it pricing} of the goods he sells : the monopolist may choose tariffs $t: Q\to \R_+$, taking into account that any agent $\tet$ buys a good $x(\tet)$ selected so as to maximize $\theta\cdot x - t(x)$; she also has the possibility of not buying anything if for every good $x$ she gets $\theta\cdot x - t(x)<0$; globally each agent buys either zero or one good (or, equivalently, she buys exactly one good but there also exists a {\it zero good} $0\in Q$ whose production cost is $0$, which gives utility $0$ to everybody, and such that the principal may not charge for it, i.e. $t(0)=0$ is prescribed). In this way the profit for the firm is 
\[\int_\Om \Big(t(x(\tet))-C(x(\tet))\Big)f(\tet)d\tet.\]
Thus, the principal's optimization problem has another optimization problem embedded: the agent's optimization problem. \\
A more standard viewpoint in {\it incentive theory}, which by the way avoids the problem of possible multiple goods $x$ maximizing the utility (net of price) of the agent, is the following : agents buy either one or zero good and, in the other hand, the firm directly chooses a family of contracts $(x(\tet),t(\tet))_{\tet\in\Omega}$ (i.e. produces $x(\tet)$ and sells it at price $t(\tet)$ to the agent of type $\tet$), and has to take into account some incentive and participative constraints. 

The first one is known as the Incentive Compatibility (IC) constraint and basically ensures that agent of type $\tet$ will really choose $x(\tet)$ at price $t(\tet)$ when she looks at the product line $(x(\tet'),t(\tet'))_{\tet' \in \Om}$. The second one is the Individual Rationality constraint (IR), modeling the fact that agents do not buy if they get a utility lower than what they would get if they bought nothing, i.e. they have a reservation utility equal to zero. We define the two constraints  as:
\begin{eqnarray}
&&\mbox{(IC) \hfill} \tet\in \mbox{argmax}\{\tet\cdot x(\tet')-t(\tet'), \; \tet' \in \Om\}  \nonumber \\
&&\mbox{(IR) \hfill} \tet\cdot x(\tet)-t(\tet) \geq 0. \nonumber
\end{eqnarray}
The classical problem then reads as
\[\sup \left\{\int_\Om \Big(t(\tet)-C(x(\tet))\Big)f(\tet)d\tet, \; (x,t) \mbox{ satisfies (IC) and (IR)}\right\}.\]
Notice that it is well-known that this problem admits a unique solution.

\bigskip
\noindent
{\bf Calculus of variations} A very powerful tool to study the above problem is the following idea : introduce (as in \cite{Arm,rc}) the potential (or {\it surplus function})
$$v(\tet):= \mbox{max} \{\tet\cdot x(\tet')-t(\tet'), \mbox{ } \tet' \in \Om\}$$ and define \[\mathcal{K_D}=\{v:\Om \rightarrow \R_+, \; v \mbox{ is convex, and } \nabla v \in \R^d_+\}.\] 
Then the problem may be re-stated in terms of $v$, in the following way :
\begin{equation}
\PD \; \; \sup_{v \in \K_D} J_D(v)= \int_\Om \Big( \tet. \nabla v(\tet) - v(\tet) - C(\nabla v(\tet)) \Big) f(\tet) d\tet. \nonumber
\end{equation}
It is a calculus of variations formulation subject to convexity constraints. This problem, after being introduced in \cite{rc}, has been studied by Carlier and Lachand-Robert in \cite{clr} as far as regularity results are concerned. Surprisingly, we will really need their  $C^1$ regularity theorem. 
Notice that this formulation is very useful to compute the solution, either explictely (see section \ref{1D}) or numerically (\cite{rc,EkeMor}). From now on we refer to the solution of $\PD$ as the {\it classical solution}, in opposition to the solution of the extended problem that we introduce in details in the next paragraph.

\bigskip
\noindent
{\bf The extended problem}
We now add a component to good qualities and to agent types. They are denoted by $y$ and $\al$, respectively, and we impose $y\geq 0$ and $\al\leq 0$. In particular, $\al$ is non-positive (we will impose that it belongs to a bounded interval of $\R_-$, $\al \in A:=[-\kappa,0]$) and it models agent's aversion to the new (one-dimensional) characteristic of the good quality.\\
In the sequel, we refer to the component $\al$ as the {\it aversion parameter} and to the corresponding one-dimensional characteristic of the good as the {\it undesirable quality} or, more simply, we call it $y$.
In the present work we look at Von Neumann Morgensten (VNM) agents and the case where the $y$ component of the good quality does not lead to any production cost.\\
In other words, we consider utilities of the type $u_{\tet,\al}(x,y,t)=\tet\cdot x+ \al y-t$, i.e. the utility that any agent of type $(\tet,\al)\in \Om\times A$ derives from buying $(x,y) \in \R^{d+1}_+$ at price $t$.\\
As before, the distribution $h$ of agent types on $\Om \times A$ is given (mathematically, it is an absolutely continuous measure on $\Om \times A$, but we will as well denote by $h$ its density). When integrating according to this density we will equivalently write $h(\tet,\al)d\tet d\al$ or $dh(\tet,\al)$. We will always tacitly assume that the $\tet-$marginal of $h$ will be the same density $f$ that we consider in the classical problem (i.e. $\int_A\! h(\tet,\al)d\al =f(\tet)$ for every $\tet\in\Omega$). 

As we mentioned before, we consider the case where the production cost (denoted in this case by $C_R$) does not depend on $y$, and we take $C_R(x,y)= C(x)$.
Introducing $w(\tet,\al):= \mbox{max} \{\tet\cdot x(\tet')+\al y(\al')-t(\tet'), \mbox{ } (\tet,\al') \in \Om\times A$ and the set 
\[\mathcal{K_R}=\{w:\Om\times A \rightarrow \R_+, \; w \mbox{ is convex, and } \nabla w \in \R^{d+1}_+\},\] we can directly jump to the variational formulation of the problem, which is
\begin{equation}
\PR \;  \sup_{v \in \K_R} J_R(w)= \int_{\Om\times A} \!\!\Big( (\tet,\al). \nabla w(\tet,\al) - w(\tet,\al) - C(\partial_{\tet} w(\tet,\al)) \Big) dh(\tet,\al). \nonumber
\end{equation}
Notice that, due to the fact that the $C-$part of the cost (the strictly convex one, which is given by the production cost) does not depend on $\partial_\alpha w$, this problem does not coincide with a multidimensional version of Mussa-Rosen's problem (i.e. with a Rochet-Choné problem; see Section 4).
We will not enter in this paper the problem of existence of a maximizer for $\PR$, which is a little bit trickier than that for $\PD$ or for Rochet-Choné, due to the lack of bounds on $\partial_\al w$. Anyway, a maximizer does exist, see for instance \cite{Car01}. 

The following - straightforward - proposition states that $\PR$ is really an extension of $\PD$.
\begin{prop}\label{relax}
If the solution $\ovw$ of $\PR$ is such that $\partial_{\al}\ovw=0$, then it coincides with the solution of $\PD$ and $\sup \PR=\sup \PD$.
In the other hand, if $h$ is a measure concentrated on $\{\al=0\}$, then $\PR=\PD$.
\end{prop}

\bigskip
In the sequel we analyse some easy features of the maximization problem $\PR$, focusing on the question whether the principal really uses positive values of the new quality $y$. The two following results show that adverse selection will ``control'' the use of $y$, even if not preventing it completely. We will interprete later these results in terms of {\it discrimination} (or {\it screening}) and {\it correlation}.
 
\begin{prop}\label{const}
If $\ovw$ is a solution of $\PR$ such that $\partial_{\al}\ovw\geq c$, for some nonnegative constant $c$, then $c=0$.
\end{prop} 
\begin{proof} 
Let $\ovw \in \K_R$ be a solution of $\PR$ such that $\partial_{\al}\ovw\geq c$. We are only concerned with perturbations of the type
\[w_{\eps}(\tet,\al)=\ovw(\tet,\al) - \eps(\al+\kappa),\]
for $0 \leq \eps\leq c$. Since $\partial_{\tet}w_{\eps}=\partial_{\tet}\ovw$ and $\partial_{\al}w_{\eps}= \partial_{\al}\ovw-\eps$, it is clear that $w_{\eps}$ is convex and increasing. The constant $\kappa$ has been inserted in order to guarantee the other admissibility constraint $w_\ve\geq 0$. To see that it is satisfied it is sufficient to check it in the corner of $\Omega\times A$ where all the coordinates are minimal. In this point we have both $\ovw=0$ and $c(\al+ \kappa)=0$, so that $w_\ve\in \K_R$. Then, we easily get
\[J_R(w_{\eps})-J_R(\ovw)= \eps \kappa \geq 0, \; \forall \eps \in [0,c].\]
This contradicts the optimality of $\ovw$ except when $c=0$ (so that $\eps=0$).
\end{proof}
\noindent

An interesting consequence of proposition \ref{const} is that the monopolist does not offer products involving nonzero $y$ when he cannot use them to discriminate (which is the case if the value of $y$ is the same for every good he sells), i.e. to screen the agents. If one restricted the offer of the monopolist adding the constraint $y=$ constant, then the solution would be the classical one.

Moreover, the same proposition proves that the minimal value of $y$ among the goods that are really present in the market must be $0$. The principal will not increase the value of $y$ for every good just for the sake of consumers' unhappiness, it will not be profitable to him!

\bigskip
There is another general case where we can prove the principal does not offer goods with undesirable qualities. Indeed, the following proposition says that he does not use the undesirable quality as soon as the distribution of  types $\tet$ and the one of aversions $\al$ are independent.
\begin{prop}
If $h(\theta,\al)=f(\tet)g(\al)$, then the classical solution is the unique minimizer of $\PR$.
\end{prop}
\begin{proof}
Assume that $\ovw$ is a solution of $\PR$. Define $v(\tet)\!=\!\!\int_{\Om \times A} \ovw(\tet,\al)g(\al)d\al$. Since $v$ does not depend on $\al$, we get $J_R(v)=J_D(v)$. On the other hand, using
\begin{gather*}
(v'(\tet))^2\leq \int \Big(\partial_\tet \ovw(\tet,\al) \Big)^2g(\al)d\al \quad\mbox{ (Jensen's inequality)},\\
-\al \partial_\al \ovw\geq 0  \quad\mbox{ (since $\al\leq 0$ and $\partial_\al \ovw\geq 0$)},\\
 \int \tet v'(\tet) f(\tet)g(\al)d\tet d\al=\int \tet \,\partial_\tet \ovw(\tet,\al) f(\tet)g(\al)d\tet d\al  \quad\mbox{ (by linearity)},
 \end{gather*}
one has:
\begin{eqnarray*}
J_R(v) \!\!\!&=&\!\!\! \int_{\Om \times A} \left(-\ud(v'(\tet))^2+\tet v'(\tet)-v(\tet) \right)f(\tet)g(\al)d\tet d\al\\
& \geq &\!\!\! \int_{\Om \times A} \!\!\!\left( -\ud (\partial_{\tet}\ovw(\tet,\al))^2\!+\!(\theta,\al)\!\cdot\!\nabla \ovw(\tet,\al)\!-\!\ovw(\tet,\al) \!\!\right)dh(\tet,\al)\!=\!J_R(\ovw).
\end{eqnarray*}
Hence, the maximization may be restricted to those convex functions that only depend on $\tet$, and, among them, the maximizer is $\bar v$. 
Also remark that the inequality  $-\al \partial_\al \ovw\geq 0$ is strict for functions truly depending on $\al$ which proves that $\bar v$ is the unique solution of $\PR$.
\end{proof}

\bigskip
In the next sections, our study consists in describing situations where of the extended problem is not the classical one. We look successively at the scalar and the multidimensional cases.
   
\section{When the principal should use the undesirable quality : the $1D$ case}\label{1D}
In this part we provide an example where the seller offers goods with various undesirable qualities $y$, hence using it for screening purpose. In particular, he uses nonzero values of $y$. 

Since Problem $\PR$ is pretty much involved in the general setting, and we are convinced that it is not easy to classify general cases where the solution of the extended problem is or is not the classical one, we first build an example in the simple case $d=1$. Notice that this means that the original problem (i.e. without $y$ and $\al$) is one-dimensional, and it actually becomes two-dimensional after adding these new variables.\\
The simplest example we can provide is the case of a quadratic production cost and a uniform distribution of type $\tet$ on segment $\Om=[\tet_{min},\tet_{max}]$. The next proposition is well-known and describes the explicit classical solution in a class of cases, which includes this one. We denote by $F$ the primitive of $f$, i.e. $F(\tet)=\int_{\tet_{min}}^\tet f(s)ds$.
\begin{prop}
Suppose that $f$ is a positive regular density on $[\tet_{min},\tet_{max}]$, with $\int_{\tet_{min}}^{\tet_{max}}f(\theta)d\theta=1$, that the production cost is quadratic, i.e. $C(x)=\frac 12 x^2$. Then, the solution $\ovv$ of problem $\PD$ is characterized by
$$\ovv(\tet_{min})=0,\quad\ovv'(\tet)=\tet-  \frac{1-F(\tet)}{f(\tet)}, $$
provided the inequalities $\tet_{min}f(\tet_{min})\geq 1$ and $f'(\theta)(1-F(\tet))\geq-2f(\tet)^2$ hold on $[\tet_{min},\tet_{max}]$. In particular, for the uniform density on $[1,2]$, these assumptions are satisfied and $\ovv(\tet)=
\frac{1}{2}\tet^2-2\tet+\frac{3}{2}$.
\end{prop}
\noindent
Although the computation of $\ovv$ is standard, we provide it for the sake of completeness.
\begin{proof}
The key point for solving $\PD$ is the integration by parts
\[\int_{\tet_{min}}^{\tet_{max}} vf=\int_{\tet_{min}}^{\tet_{max}}v'(1-F),\]
where we used $F(\tet_{max})=1$ and $v(\tet_{min})=0$ (we may assume $\inf v =0$ and, $v$ being increasing, the minimum si realized at $\tet_{min}$).
We set $u=v'$ and we have:
\[J_D(v)=\tilde{J}(u)=\int_{\tet_{min}}^{\tet_{max}} \Big(\ud u(\tet)^2-\tet u(\tet)\Big)f(\tet) +u(\tet)\Big(1-F(\tet)\Big) d\tet .\]
Then we easily get that the minimal value would be realized by choosing:
\[u(\tet)=\tet -\frac{1-F(\tet)}{f(\tet)}, \; \tet \; a.e.,\]
provided this choice is possible. This means that $u(\tet)$ must be the derivative of a function $v$ satisfying the constraints. Since $v$ must be increasing and convex, then $u$ must be positive and increasing. It is sufficient to compute its derivative, to require that it is positive on $[\tet_{min},\tet_{max}]$, and to impose $u(\tet_{min})\geq0$. This last condition gives $\tet_{min}-\frac{1}{f(\tet_{min})}\geq 0$, and the condition on $u'$ gives
$$0\leq u'(\tet)=1-\frac{-f(\tet)^2+(1-F(\tet))f'(\tet)}{ f(\tet)^2}=2+\frac{(1-F(\tet))f'(\tet)}{ f(\tet)^2},$$
which is known as the {\it hazard condition}. 
The assumptions on $F$ exactly ensure that these conditions are satisfied.
For the particular case of the uniform density on $[1,2]$, just replace $F(\tet)=\tet-1$ and $f(\tet)=1$ and deduce from this equation that $\ovv'(\tet)=2(\tet-1)$. We then obtain the expression of $\ovv$ above and finally check that $\ovv \in \K_D$. 
\end{proof}
\noindent
Note that it is a simple case without bunch (or equivalently perfect screening i.e. the monopolist offers different contracts to different types).

\bigskip
Let us now look at the extended problem $\PR$ for the quadratic production cost. Our goal is to find some cases where the solution of $\PR$ is not the classical solution $\ovv$ of $\PD$.\\ 
The strategy starts by introducing, for all $\eps>0$, the perturbations \[w_{\eps}(\tet,\al)=\ovv(\tet+\eps l(\al)), \; \forall (\tet,\al) \in \Om \times A,\]
where $ l(\al)=(\al+a)_+=\max(\al+a,0)$ and $a \in (-\kappa,0)$ is given. If necessary, one can extend $\bar v$ outside $\Omega$ by affine extensions before defining $w_\ve$. Notice that $w_{0}=\ovv$ and that, due to the convexity and monotonicity properties of $l$, $w_\ve$ is convex and increasing as well.

\begin{figure}[!h]
   \centering
   \input{expl1D}
   \caption{Structure of the $1D$ example}\label{R1R2}
\end{figure}
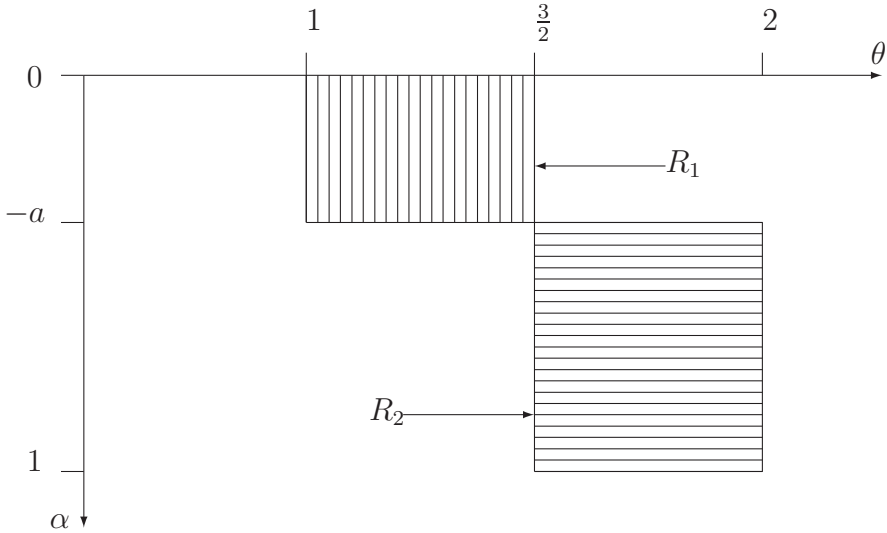

\begin{prop}\label{prop_1D}
Let $\ovv$ be the classical solution of $\PD$ for a certain density of types $f$. Assume that $\ovv$ is $C^2$ and that there exists $B \subset \Om$ such that $\int_B f(\tet)d\tet>0$ and $\forall \tet \in B, \; \ovv''(\tet)\left(\tet- \ovv'(\tet)\right)-\ovv'(\tet)>0$.  Then there exists a joint distribution $h$ on $\Om \times A$, with $\tet-$marginal equal to $f$, such that $J_R(w_\ve)>J_R(\bar v)$ for small $\ve>0$. In particular $\sup \PR < \sup \PD$ and the solution of $\PR$ is not $w(\tet,\al)=\ovv (\tet)$.\\
An explicit example is obtained if one takes  $\Omega=[1,2],\;\kappa=1, \; a\in (0,1), \; f=1$ and $h(\tet,\al)=\frac{1}{a}\mathrm{1}_{(\tet,\al) \in R_1}+\frac{1}{1-a}\mathrm{1}_{(\tet,\al) \in R_2}$, where $R_1$ and $R_2$ are the two regions of Fig.  \ref{R1R2}.
\end{prop}
\noindent
Since $J_R(\ovv)=J_D(\ovv)=\sup \PD$ and using proposition \ref{relax} we immediately deduce the next corollary.
\begin{coro}
For some joint distributions of types and aversions, the principal offers goods with nonzero undesirable quality $y$.
\end{coro}
\noindent

The previous result states that the principal takes advantage of the costless and undesirable quality to discriminate. 
 Indeed, screening considerations and the possibility to ease the incentive constraints outperform the desire to provide every agent with fully desirable goods ($y=0$). \\Before proving proposition \ref{prop_1D}, we would like to comment a bit such particular case. The coupling between the classical type $\tet$ and the aversion characteristic $\al$ is concentrated on the two regions $R_1$ (low types have low aversion parameters) and $R_2$ (high types have  higher aversion parameters) represented in Fig. \ref{R1R2}. From the proof we will infer that in the present example the principal derives profit from offering undesirable goods to those agents having both low types $\tet$ and low aversion $\al$ (i.e. in $R_1$). Hence, in the present case, a correlation aspect appears: the example we provide has the same structure than the one obtained in \cite{mr_1} and \cite{BolDew} in a discrete framework, where optimal contracts are random when the correlation between $-\al$ and $\tet$ is positive. \\
We now prove proposition  \ref{prop_1D}.

\begin{proof}
First we compute the derivatives of $w_{\eps}$ : $\partial_{\tet}w_{\eps}=\ovv'$ and $\partial_{\al}w_{\eps}=\eps l'\ovv'$, where $l$ is differentiable (i.e. almost everywhere). We already noticed the admissibility of $w_{\eps}$, which follows from the convexity and monotonicity of $l$ and $\ovv$.
Now, define $j(\eps):=J_R(w_{\eps})$ and compute the first variation
\begin{eqnarray}\label{int_j}
j'(0)\!\!=\!\!\int_{\Om \times A} \left[  \Big(\ovv''(\tet)\left(\tet- \ovv'(\tet)\right)-\ovv'(\tet)\Big)l(\al)+\al \ovv'(\tet)l'(\al)\right] dh(\tet,\al). 
\end{eqnarray}
Assume that there exists a set $B$ satisfying the same assumptions as in the proposition. Then we can choose a density $h$ of the following form : $h=\frac{1}{\tilde a}f(\tet)\mathrm{1}_{[-\tilde a,0]}+\frac{1}{\kappa-a}f(\tet)\mathrm{1}_{[-\kappa,-a]}(\al)$, for some small value $\tilde{a} \in (0,a]$. In this case it is clear that the integral in \pref{int_j} will be immediately restricted to the set $B \times (-\tilde{a},0]$, where the first term is positive. Since the goal is to get $j'(0)>0$, we only need to choose $\tilde a$ small enough. In this case, the term involving $\al \bar v'(\tet) l'(\al)$ will be negligible (since the density $h$ tends to concentrate around the set $\{\al=0\}$), and  
$$\lim_{\tilde a\to 0} \!\!\int \!\!  \Big(\ovv''(\tet)\left(\tet- \ovv'(\tet)\right)-\ovv'(\tet)\Big)l(\al)dh\!= \!l(0)\!\!\int\!\! (\ovv''(\tet)\left(\tet- \ovv'(\tet)\right)f(\tet)d\tet \!>\!0.$$
Thus, we finally get  $j'(0)>0$.

\bigskip
\noindent
To prove the second statement  of the proposition, we simply compute the variation replacing $\ovv$ (who is of course $C^2$ in the uniform case) and $h$ by their expression 
$$
j'(0)=\int_{\Om \times A} \Big( (3-2\tet)l(\al)+\al(\tet-1)l'(\al)\Big)dh(\tet,\al)= \frac{a^2}{16} > 0.
$$
Notice that in this particular case, $B=[1,\frac{3}{2}]$ and $\tilde{a}=a$.
\end{proof}

\begin{rem}
We emphasize that, in the case where the production cost in the extended problem has the form $C_R(x,y)= C(x)+\lambda y$ for some $\lambda>0$ 
one can anyway build an explicit example from the uniform density $f$. More precisely, for all $\lambda>0$, there exists $M=M_{\lambda}>0$ and a density $h$ on $[1,2]\times [-M,0]$ providing a non-classical solution. This density may be chosen of the usual form, taking $0<a<M$, with $\lambda<\frac12 a$, $h(\tet,\al)=\frac{1}{a}\mathrm{1}_{(\tet,\al) \in B\times [-a,0]}+\frac{1}{M-a}\mathrm{1}_{(\tet,\al) \in \Om \backslash B \times [-M,-a]}$. Notice anyway that for a fixed set of aversion parameter $A$ this is possible for some $\lambda$ only.
\end{rem}

\begin{rem}
Note that the assumption that $\ovv$ should be $C^2$ is not really needed to get such a result. Indeed we could suppose that $\ovv$ is only $C^1$ and apply a similar technique as the one we use for proving proposition \ref{prop_2D} in the next section.
\end{rem}

\section{The multidimensional case}\label{2D}

In this section we look at the multidimensional case, which has been mainly studied in \cite{Arm,rc}. Since less information on the solution in this case is known than in the 1D one, we will try to use the minimal regularity and the minimal structure results on the optimal $\bar v$ that have been investigated. 

\begin{figure}[!h]
   \centering
   \input{expl2D}
   \caption{The shape of the classical solution for the exponential distribution}\label{ex2D}
\end{figure}
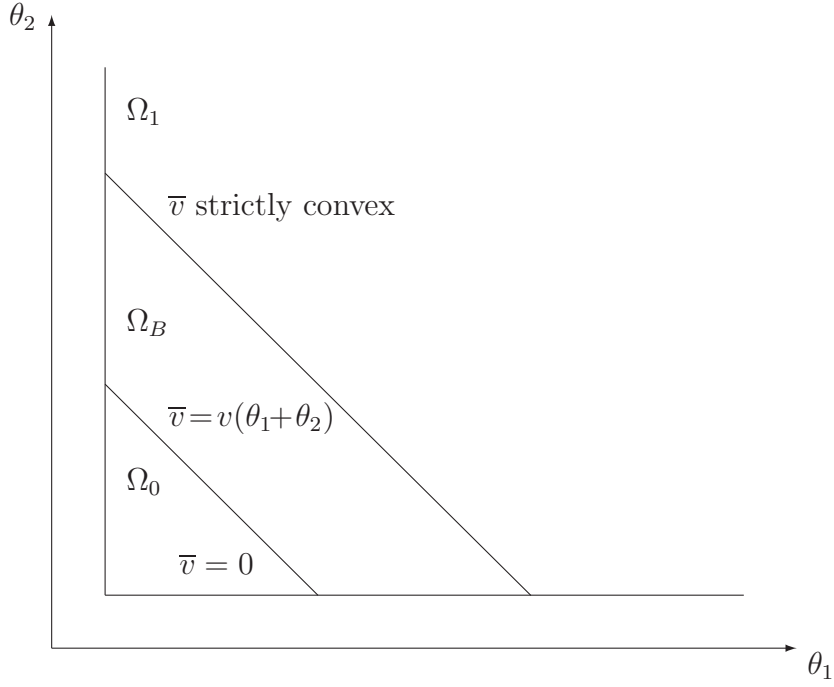

\begin{prop}\label{prop_2D}
Suppose that $f$ is a density on a convex domain $\Omega\subset\R^d$ such that the optimal function $\bar v$ in the Rochet-Choné problem $\PD$ is a $C^1$ convex function, with $\Omega_0:=\{\bar v=0\}$ having non-empty interior. Then there exists a density $h$ on $\Omega\times A$ such that the solution of $\PR$ is not $w(\theta,\al)=\bar v(\tet)$.
\end{prop}
\noindent
For the sake of this proposition, we need two things: $C^1$ regularity of $\bar v$ and the fact that there exists a flat zone $\{\bar v=0\}$.\\
General $C^1$ regularity results are proven in \cite{clr} for the solutions of a wide class of variational problems under convexity constraints, including the one of Rochet \& Choné. As far as the structure of $\ovv$ is concerned, \cite{rc} proves that if $d=2$, $\bar v$ has the shape described in Fig.\ref{ex2D} when $f$ is an exponential density, and check numerically the same behavior for the uniform density. For numerical computations showing this kind of behavior, we also refer to \cite{EkeMor}. Anyway the crucial feature of the classical solution we need is that the region $\Omega_0:=\{\bar v=0\}$ has non-empty interior. This statement means that the seller excludes a set of consumers. Actually, it has been proven by Armstrong in \cite{Arm} that for any dimension $d \geq 2$, the set $\Omega_0$ has positive measure as soon as $\Omega$ is strictly convex; yet, the same proof easily generalizes to other classes domains, for instance to any straight rectangle (i.e. a rectangle whose edges are parallel to the coordinate axes).  We thus easily derive the following corollary.
\begin{coro}\label{coro_D}
For every regular density $f$ of types on a domain $\Omega$ which is either strictly convex or a straight rectangle in $\R^d$, for $d\geq 2$, there exists a density $h$ on $\Omega\times A$ such that the optimal contract contains undesirable goods.
\end{coro}

We will concentrate the proof of proposition \ref{prop_2D} in the case of a straight rectangle in 2D,  which is the one that has been most investigated so far, but the reader may easily notice that the procedure may be extended to higher dimensions and to different shapes of $\Omega$, since we only perform local variations.

\begin{proof}
The strategy for looking at the multidimensional case is similar to the 1D one. We suppose that the solution is a function $\ovv$ depending on $\theta$ only, and we produce a function $w_\ve$ through $w_\ve(\theta,\al)=\bar v(\tet+\ve \ev l(\al))$, where $l$ is the usual function $l(\al)=(\al+a)_+$ and $\ev$ is a given unit vector with positive components. In the following we choose $\ev=(\frac{1}{\sqrt{2}},\frac{1}{\sqrt{2}}).$ The function $w_\ve$ is convex and increasing with respect to each variable, provided that $\ovv$ is convex and increasing. We notice that 
$$\nabla_\tet w_\ve(\tet,\al)=\nabla\bar v(\tet+\ve \ev l(\al)),\; \partial_\al w_\ve= \ve \nabla\bar v(\tet+\ve \ev l(\al))\cdot \ev l'(\al),\;\mbox{ for a.e. }(\tet,\al).$$
Hence we can write $J_R(w_\ve)$ and we have
\begin{eqnarray*}
J_R(w_\ve)\!\!\!\!&=&\!\!\!\int \!\!\left(\!\theta\!\!\cdot\!\! \nabla\bar v(z_{\tet,\al})+\ve \alpha l'(\al) \nabla\bar v(z_{\tet,\al})\!\!\cdot\!\ev -\ovv(z_{\tet,\al})-\!\frac{|\nabla\bar v|^2}{2}(z_{\tet,\al})\!\!\right)\!dh(\tet,\al)\\
&=&\!\!\!\int \left(-k(z_{\tet,\al})-\ve(l(\al)- \alpha l'(\al)) \nabla\bar v(z_{\tet,\al})\cdot\ev \right)dh(\theta,\al)\end{eqnarray*}
where $z_{\tet,\al}$ denotes $\tet+\ve \ev l(\al)$ and  $k(z):=\frac 12 |\nabla\bar v(z)|^2-z\cdot\nabla\bar v(z)+\bar v(z)$.

\smallskip
\noindent{\bf Simplifying assumptions on $h$}\\
Since our goal is just to provide sufficient conditions to show an example where $\bar v$ is not optimal (and $w_\ve$ gives a better result than $\bar v$), we can make some simplifying assumptions: we suppose that $h(\tet,\al)=I_K(\tet) I_S(\al)$ for those pairs $(\tet,\al)$ such that $l(\al)\neq 0$, where $K\subset \Omega$ is a non-negligible given subset of $\Omega$ and $S\subset [-a/2,0]$ is an interval where $l$ is positive. This is not so restrictive in terms of the density $f$, provided $f$ is bounded from below by a positive constant on $K$: it is actually always possible to express $f$ as the $\tet-$marginal of a density $h$ of this kind, if we adjust the rest of the density with the behavior of $h$ on $\{l=0\}$. This reduces the integral above to 
$$-\int_S d\al \int_K \left(k(z_{\tet,\al})+\ve(l(\al)- \alpha l'(\al)) \nabla\bar v(z_{\tet,\al})\cdot\ev \right)d\tet.$$

\smallskip
\noindent{\bf Simplifying assumptions on $K$}\\
Now suppose that $K$ is actually a rectangle with two sides which are parallel to $\ev$ and two orthogonal sides denoted by $\partial K_+$ and $\partial K_-$ (see Fig. \ref{struct2D}). If we denote by $\,d_1$ the one-dimensional integration on $\partial K_+$ and $\partial K_-$, it is not difficult to check that 
\begin{eqnarray}
\frac{d}{d\ve}\Big|_{\ve=0} \int_K k(\tet+\ve \ev l(\al))d\tet&=&\frac{d}{d\ve}\Big|_{\ve=0}  \int_{K+\eps \ev l(\al)}\!\!\!\!\!k(\tet)d\tet \nonumber \\
&=&l(\al) \Big( \int_{\partial K_+} \!\! \!\!k(\tet) \,d_1\tet-\int_{\partial K_-} \!\!\!\!k(\tet) \,d_1\tet \Big)\nonumber
\end{eqnarray}
provided that $k$ is continuous (which is the case if $\ovv\in C^1$). Moreover,  we have
$$\int_K \nabla\bar v\cdot \ev\,d\tet=\int_{\partial K_+} \bar v \,d_1\tet-\int_{\partial K_-} \bar v \,d_1\tet.$$
Now, suppose that $\bar v=\nabla\bar v =0$ on $\partial K_-$. This simplifies the computation and yields to 
$$\frac{d}{d\ve}\Big|_{\ve=0} J_R(w_\ve)\!=\!\int_S  l(\al)d\al\left(\int_{\partial K_+}\!\! k \,d_1\tet\!\right)+\int_S \!\left(l(\al)-\alpha l'(\al)\right)d\al \left(\int_{\partial K_+} \!\!\bar v \,d_1\tet\!\right)\!.$$
Since $\ovv$ is nonnegative and $l(\al)-\alpha l'(\al)=(\al+a)_+-\al=a\leq 2l(\al)$ on $S$, we can estimate the last term through \[\int_S \left(l(\al)-\alpha l'(\al)\right)d\al \left(\int_{\partial K_+} \bar v \,d_1\tet\right)\leq 2\int_S l(\al)d\al  \left(\int_{\partial K_+} \bar v \,d_1\tet\right),\] and get
$$\frac{d}{d\ve}\Big|_{\ve=0} J_R(w_\ve)\geq - \int_S  l(\al)d\al\left(\int_{\partial K_+}\!( k+2\bar v) \,d_1\tet\right).$$

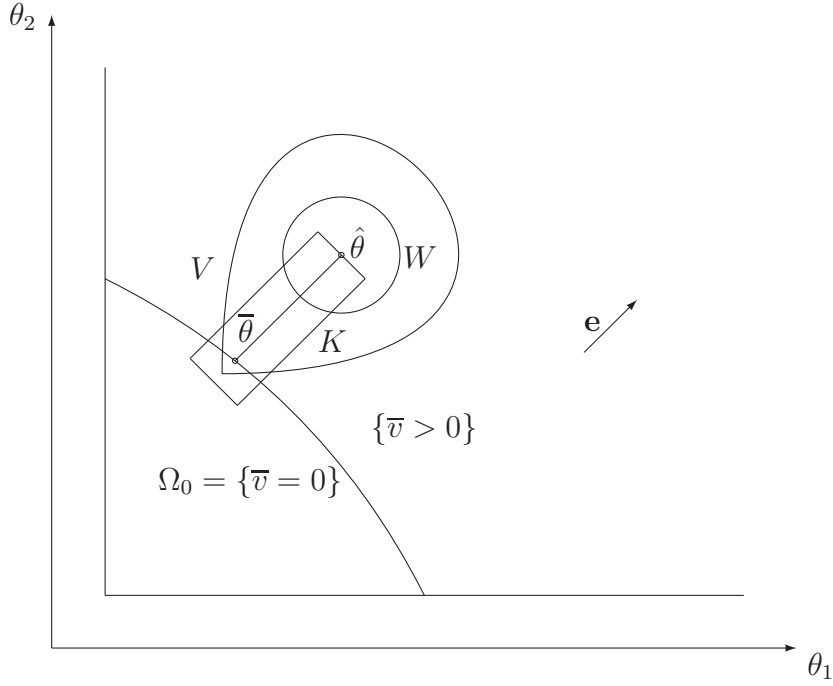
\begin{figure}[!h]
\begin{center}
\setlength{\unitlength}{1.4cm}
\begin{picture}(7,5.6)(0,0.5)
   \put(0.5,0.5){\vector(1,0){7}}
 
   \put(5.5,3.3){\vector(1,1){0.5}}
      \put(5.5,3.5){$\ev$}
   \put(1,1){\line(1,0){6}} 
   \put(7.6,0.25){$\tet_1$}
    \put(0.5,0.5){\vector(0,1){6}}
    \put(1,1){\line(0,1){5}} 
    \put(0.1,6.4){$\tet_2$}
   \qbezier(1,4)(3,3)(4,1)
\cbezier(2.1,3.1)(2.1,8.2)(7.1,3.1)(2.1,3.1)
\put(3.222,4.222){\circle{1.1}}
    \put(1.5,2){$\Omega_0=\{\ovv=0\}$} 	
     \put(3.5,2.5){$\{\ovv>0\}$} 	
    \put(2.25,3.4){$\bar\tet$} 	
     \put(2.222,3.222){\line(1,1){1}} 
       \put(3.3,4.2){$\hat\tet$} 
         \put(3.8,4.1){$W$} 
 \put(1.8,4){$V$} 
\put(2.222,3.222){\circle{.05}}
\put(3.222,4.222){\circle{.05}}
              \put(3,3.3){$K$} 
         \put(3.444,4){\line(-1,1){0.444}} 
 \put(2.244,2.8){\line(-1,1){0.444}} 
       \put(2.244,2.8){\line(1,1){1.2}} 
        \put(1.8,3.244){\line(1,1){1.2}} 
\end{picture}
\end{center}
\caption{Structure of the $2D$ example}\label{struct2D}
\end{figure}

\smallskip
\noindent{\bf Choice of a point $\bar\tet$}\\
It is now sufficient to prove that one can choose $K$ so that $k+2\bar v= \frac 12 |\nabla\bar v(\tet)|^2-\tet\cdot\nabla\bar v(\tet)+3\ovv(\tet)$ is negative on $\partial K_+$ in order to have a positive derivative of $J_R(w_\ve)$.

Here the key point is to suppose that $\ovv$ has a flat part $\Omega_0$ with non-empty interior. Thanks to the convexity of $\bar v$, this set is convex and we suppose that it includes a neighborhood of $0$. 
Let $\bar\theta$ be a point on the boundary of $\Omega_0$ which is in the interior of $\Omega$. Then, thanks to the monotonicity properties of $\ovv$, the vector $\ev$ is an outwards vector for the set $\Omega_0$, in the sense that $\bar v(\bar\tet+ s\ev)>0$ for $s>0$.  The vector $\frac 12 \nabla\bar v(\bar\tet)-\bar\tet$ will have negative coordinates, and the same is true in a neighborhood $V\subset\Omega$ of $\bar\tet$. If $-c$ is a constant bounding from above the components of $\frac 12 \nabla\bar v(\tet)-\tet$ for $\tet\in V$,  we can write
$$\left(\frac12 \nabla\bar v(\tet)-\tet\right)\cdot\nabla\bar v(\tet)\leq -c\ev\cdot \nabla\bar v(\tet),$$
and the term $\ev\cdot \nabla\bar v(\tet)$ is the derivative of $\ovv$ in the direction $\ev$.

\smallskip
\noindent{\bf Choice of $\partial K_+$ so that $k+2\bar v<0$}\\
Now, let us consider the values of $k+2\bar v$ on a line segment stemming from $\bar\tet$ in the direction $\ev$, i.e. we consider $(k+2\bar v)(\bar\tet+s\ev)$ for $s>0$ such that $\bar\tet+s\ev\in V$. Let us suppose that this quantity stays always positive : this means that we have the inequality 
\begin{eqnarray*}
-c\frac{d}{dt} \!\!\Big(\!\bar v(\bar\tet\!+\!s\ev)+2\bar v(\bar\tet\!+\!s\ev)\!\Big)\!\!\!\!\!&\geq &\!\!\!\!\left(\frac12 \nabla\bar v(\bar\tet+s\ev)-(\bar\tet\!+\!s\ev)\!\!\right)\!\!\cdot\!\!\nabla\bar v(\bar\tet\!+\!s\ev)\!\!+\!\!2\bar v(\bar\tet\!+\!s\ev)\\
&\geq& 0\;,
\end{eqnarray*}
which implies that
$$\frac{d}{dt} \bar v(\bar\tet+s\ev)\leq \frac{2}{c} \bar v(\bar\tet+s\ev).$$
By Gronwall Lemma, this inequality implies that $\bar v(\bar\tet+s\ev)\leq e^{\frac 2c s} \bar v(\bar\tet)=0$, which is in contradiction with the fact that $\bar\tet$ is on the boundary of $\Omega_0$, i.e. that $\bar v(\bar\tet+ s\ev)>0$.

This implies that there exists a point $\hat\tet$ of the form $\bar\tet+ s\ev$ in the neigborhood $V$ where $k+2\bar v$ is negative. The same will be true for all the points in a neighborhood $W$ of $\hat\tet$. In this way, if one chooses a rectangle $K$ such that one of the two sides orthogonal to $\ev$ is completely included in $\Omega_0$ and the other is included in $W$, the derivative of $J_R(w_\ve)$ is positive.
\end{proof}
Notice that in this proof, we use again profitable (from the seller viewpoint) contracts where goods with positive $y$ are sold to rather low type agents (i.e. the lowest type agents among those who enter the market). Actually, as one can see from Fig. \ref{struct2D}, the agents who are concerned are those who are close to the region $\Omega_0$, i.e. to the consumers who stay out of the market.

\section{Summary and conclusions}\label{sconc}

In this paper, we have been concerned with providing partial progress on the possible profitability of random (or undesirable) contracts, for a monopolist facing risk averse (or averse) buyers with multidimensional types. In comparison with the existing literature which is focused on the random interpretation,  we stressed the fact that, in this simplified model, randomness is just an instance of undesirable qualities. This is why we feel free to switch from one terminology to the other.

\bigskip
We first defined an extended version of the multidimensional screening problem mainly studied by Rochet \& Chon\'e. We discussed several interpretations of the addition of an undesirable quality and we highlighted situations where the seller either can or cannot improve his profit by offering undesirable goods (or lotteries). 

\bigskip
The optimal policy is not random in the case where agents type $\tet$ and aversion $\al$ are independent. The correlation between the distributions of agents type and aversion seems to be at the heart of the answer. We built profitable examples involving the undesirable quality (risk, say) when the  correlation is of the following form: low-type consumers are less risk averse. This can be seen both in the 1D and in the multidimensional examples we provided, where we exploited the presence of poorly risk averse agents close to the ``out-of-the-market'' zone $\Omega_0$. What we get is consistent with the results presented in \cite{BolDew,mr_1}. Yet, the novelty of this paper is to analyze these same questions with a continuum of agents, in a Rochet \& Choné framework.

\bigskip
The non profitability conclusion is also true when the monopolist is added a constraint forcing him to offer constant risk, i.e. removing the possibility to screen the agents (or discriminate). We saw that, if the principal really uses random contracts (and he only does use them for pure screening purposes), then the randomness he proposes (i.e. the undesirable characteristic $y$) will take different values, including $0$. \\
Anyway, by providing examples in both the scalar and the multidimensional cases, we emphasized the role of the undesirable quality to screen the agents, and we have seen that the screening considerations sometimes outperform aversion.

\subsection*{Acknowledgements}
The authors wish to thank G. Carlier for suggesting the topic of this article as well as for very helpful discussions and I. Ekeland for his active interest in the problem and interesting suggestions. The second author acknowledges the support of the ANR project ``Evolution and variational methods applied to economics and finance''.

\end{document}

%% file: expl1D.tex
%
%
%

\setlength{\unitlength}{1.5cm}
\begin{picture}(7,5)(0,-3.5)
   \put(0,1){\vector(1,0){7}}
   \put(1.95,1){\line(0,1){0.2}} 
   \put(3.95,1){\line(0,1){0.2}}
   \put(5.95,1){\line(0,1){0.2}}
   \put(1.95,1.4){1} 
   \put(3.95,1.4){$\frac{3}{2}$}
   \put(5.95,1.4){2}
   \put(6.9,1.1){$\tet$}
    \put(0,1){\vector(0,-1){4}}
    \put(-0.2,1){\line(1,0){0.2}} 
    \put(-0.2,-.3){\line(1,0){0.2}}
    \put(-0.2,-2.5){\line(1,0){0.2}}
    \put(-0.5,0.9){$0$} 
    \put(-0.7,-.3){$-a$}
    \put(-0.5,-2.5){1}
   \put(-0.3,-3){$\al$}
   \put(1.95,-.3){\line(1,0){4}}
   \put(3.95,-2.5){\line(1,0){2}}
   \put(3.95,1){\line(0,-1){3.5}}
   \put(1.95,1){\line(0,-1){1.3}}
   \put(5.95,-.3){\line(0,-1){2.2}}
   \multiput(1.95,1)(0.1,0){21}{\line(0,-1){1.3}} 
   \multiput(3.95,-.3)(0,-0.1){23}{\line(1,0){2}} 
   \put(5.1,0.2){\vector(-1,0){1.15}}
   \put(2.8,-2){\vector(1,0){1.15}}
   \put(2.5,-2.05){$R_2$}
   \put(5.1,0.15){$R_1$}
\end{picture}


%% file: expl2D.tex
%
%
%

\setlength{\unitlength}{1.4cm}
\begin{picture}(7,7)(0,0)
   \put(0.5,0.5){\vector(1,0){7}}
   \put(1,1){\line(1,0){6}} 
   \put(7.6,0.25){$\tet_1$}
    \put(0.5,0.5){\vector(0,1){6}}
    \put(1,1){\line(0,1){5}} 
    \put(0.1,6.4){$\tet_2$}
    \put(1,3){\line(1,-1){2}}   
    \put(1,5){\line(1,-1){4}} 	
    \put(1.2,2){$\Omega_0$} 	
    \put(1.7,1.2){$\bar v =0$} 	 
    \put(1.2,3.5){$\Omega_B$} 	
        \put(1.6,2.6){$\bar v \!= \!v(\tet_1\!\!+\!\tet_2)$} 	
    \put(1.2,5.5){$\Omega_1$} 	
    \put(1.6,4.6){$\bar v$ strictly convex} 
\end{picture}
